\documentclass[fleqn,12pt]{article}

\usepackage{amsmath, amscd, amsfonts, amssymb, graphicx, color}
\usepackage{hyperref}
\usepackage{cite}
\usepackage{amsmath}
\usepackage{amsfonts}
\usepackage{epsfig}
\usepackage{amssymb}
\usepackage{amsthm}
\usepackage{epstopdf}
\usepackage{footnote}
\usepackage{newlfont}
\usepackage{color}

\begin{document}
%%%%%%%%%%%%%%%%%%%%%%%%%%%%%%%%%%%%%%%%%%%%%%%%%%%%%%%%

 \title{Digital image watermarking using normal matrices}

 \author{E. Kokabifar\thanks{Department of Mathematics,
 Faculty of Science, Yazd University, Yazd, Iran
 (e.kokabifar@stu.yazd.ac.ir, loghmani@yazd.ac.ir).},\,
 G.B. Loghmani\footnotemark[1]\,\,
 and  A. Latif\thanks{Department of Electrical and Computer Engineering,
  Yazd University, Yazd, Iran (alatif@yazd.ac.ir).}}
\maketitle

\vspace{-6mm}

\begin{abstract}
This paper presents techniques for digital image watermarking based on eigenvalue decomposition of normal matrices. The introduced methods are convenient and self-explanatory, achieve satisfactory results, as well as require less and easy computations compared to some current methods. Through the proposed methods, host images and watermarks are transformed to the space of normal matrices, and the properties of spectral decompositions are dealt with to obtain watermarked images. Watermark extraction is carried out via a procedure similar to embedding. Experimental results are provided to illustrate the reliability and robustness of the methods. 
\end{abstract}


%%%%%%%%%%%%%%%%%%%%%%%%%%%%%%%%%%%%%%%%%%%%%%%%%%%%%%%%%%%%%%%%%%%%%%%
\section{Introduction}\label{intro}

Nowadays, everyone has an easy access to digital images via the Internet and other means of communication. In addition, powerful software and hardware tools are available to edit any kind of digital multimedia. So, digital images may become a problematic issue if it have significant authentications and information. To cope with this problem, \textit{digital image watermarking} can be employed. Digital image watermarking is a procedure of embedding additional information called \textit{watermark} into an image, called \textit{host image}, which preserves the perceptual qualities of the original host image. When used normally, a watermark, carrying important information, is hidden and can be visible only by pursuing some particular processes. The watermark should be detectable or extractable from the host image by owner identification or integrity verification of the watermarked image \cite{latif}.

Watermarking has a long list of applications in such areas as protecting copyright ownership of digital media, verifying the authenticity of a bill or the trademark of a paper manufacturer, making the information imperceptible, keeping the information secret, and so on. Digital image watermarking has received a great deal of attention as much as to be practiced in many other important applications. References for theory and practice of this method are  \cite{latif0,survey,review} and \cite{lin}, to name but a few.

For the subject to be fully grasped, presenting certain basic definitions seems necessary. A \textit{gray scale image} can be demonstrated as a matrix the magnitude of each of whose components is called \textit{gray level} or \textit{intensity} of the image at that point. A given image is called a digital image if all of its elements are finite and discrete quantities; in this case, these elements are called \textit{pixels}. \textit{Digital Image Processing} (DIP) concerns using digital computers to preform processing on digital images. DIP receives a digital image as an input, processes efficient algorithms on it by a computer and gives an image as an output. So, digital image watermarking can be construed as a DIP.
 
The theory of linear algebra provides helpful instruments for digital image watermarking as well as a variety of other fields of digital image processing. One of the most useful tools for digital image watermarking is employing Singular Value Decomposition (SVD) \cite{chung,chang1}.

the matrix $A \in \mathbb{C}^{m\times n}$ can be written in the form of $USV^{*}$, where $U\in\mathbb{C}^{m\times m}$ and $V\in\mathbb{C}^{n\times n}$ are unitary matrices, i.e., $U^*U=I_m, V^*V=I_n$, where $^*$ denotes complex conjugate transpose. The matrix $S$ is an $m\times n$ diagonal matrix such that its nonnegative entries are ordered in a non-increasing order (see for example, Theorem 7.3.5 of \cite{horn}). Motivated by the influences and traits of singular values of $A$ in digital image watermarking, and regarding the point that the singular values of $A$ are the positive square roots of the eigenvalues of matrices $AA^*$ and $A^*A$, this paper concerns itself with the eigenvalue of the normal matrices $A+A^*$ and $A-A^*$ in order to obtain techniques for digital image watermarking which are more effective, yield better results and need fewer computations.

In the next section, some definitions and concepts are briefly presented about normal matrices. Section \ref{water} consists of two subsections in which the proposed image watermarking methods are explained. In Section \ref{exp}, the validity rates of the presented watermarking schemes are investigated, and their efficiencies are compared by experimental results. In addition, the robustness of watermarking is evaluated against some common attacks and criticisms.

\section{Normal matrices} \label{normal}
In this section, definitions and some properties of normal matrices are reviewed. See \cite{normal1,normal2} and references therein as the suggested sources on a list of conditions on normal matrices. The proposed method based on these presented properties is explained in the next section.

A matrix $M\in\mathbb{C}^{n \times n}$ is called \textit{normal} if $M^*M=MM^*$. Assume that $M$ is an $n$-square normal matrix, then there exists an orthonormal basis of $\mathbb{C}^{n\times n}$ consisting of eigenvectors of $M$, and $M$ is unitarily diagonalizable. Let the scalars $\lambda_1,\ldots, \lambda_n$, counted according to multiplicity, be  eigenvalues of the normal matrix $M$ and let $u_1, \ldots, u_n$ be its corresponding orthonormal eigenvectors. Then, the matrix $M$ can be factored as
\begin{equation*}
M=U\Lambda U^*,\qquad \Lambda=\mbox{diag}(\lambda_1,\ldots,\lambda_n),\qquad U=[u_1,\ldots,u_n],
\end{equation*}
where the matrix $U$ satisfies $UU^*=I_n$. Without loss of generality, assume that eigenvalues are ordered in a non-ascending order of magnitude, i.e., $\left| {\lambda _1 } \right| \ge \left| {\lambda _2 } \right| \ldots  \ge \left| {\lambda _n } \right|$. 

Note that, if all the elements of the matrix $M$ are real, then $M^*=M^T$, where $M^T$ denotes the transpose of the matrix $M$.  A square matrix $M$ is called \textit{symmetric} if $M=M^T$ and called \textit{skew-symmetric} if $M=-M^T$. It is easy to see that symmetric and skew-symmetric matrices are normal. Also, all the eigenvalues of a real symmetric matrix are real, and all the eigenvalues of a real skew-symmetric matrix are  purely imaginary.
A general square matrix $M$ satisfies $M=B+C$, where the symmetric matrix $B=\dfrac{M+M^T}{2}$ is called \textit{symmetric part} of $M$ and, analogously, the skew-symmetric matrix $C=\dfrac{M-M^T}{2}$ is called \textit{skew-symmetric} part of $M$. Consequently, every square matrix can be written as the sum of two normal matrices: a symmetric matrix and a skew-symmetric one. This point is specially used in the proposed watermarking schemes.

\section{Digital image watermarking method}\label{water}
This section includes two subsections in which we present methods for digital image watermarking using normal matrices. To do this, we transform both the matrix representing the host image and the one representing the watermark into the space of normal matrices. Next, we benefit from the properties of its eigenvalue decompositions and embed the watermark into the host image. Finally, returning to our original space, we can construct the watermarked image. Also, we may extract the watermark by following almost the same processes.

Let $X$ and $W$ be two square matrices representing the host image and the watermark image respectively. Two distinct methods are considered. First, we deal with the symmetric parts of $X$ and $W$ to obtain a digital image watermarking scheme. Next, we explain another technique using both symmetric and skew-symmetric parts of the matrices $X$ and $W$. It is noticeable that finding the eigenvalues and eigenvectors of a matrix needs less computation than finding its singular values and singular vectors. Furthermore, the eigenvalues and eigenvectors of a normal (especially symmetric or skew-symmetric) matrix can be calculated by explicit formulas and, therefore, may yet again require less computation \cite{sym1,sym2,sym4}.

\subsection{Digital image watermarking method using symmetric parts of $X$ and $W$}\label{sb1}
In this subsection, we describe a digital image watermarking scheme based on the eigenvalue decomposition of the symmetric parts of the matrices $X$ and $W$. This procedure can be carried out similarly for the case of the skew-symmetric parts of the two matrices $X$ and $W$. Notice that, for the same host image and watermark, results obtained by these two techniques (using symmetric or skew-symmetric parts) may be different.

Assume that $X$ is an $n \times n$ matrix and $W$ is an $m \times m$. Let $B_X$ and $B_W$ be the symmetric parts of the two matrices $X$ and $W$, respectively. Therefore, the normal matrices $B_X$ and $B_W$ can be factored as follows
\[
B_X  = U_{B_X } \Lambda _{B_X } U_{B_X }^* ,\qquad\Lambda _{B_X }  = \mbox{diag}(\lambda _{B_X ,1} , \ldots ,\lambda _{B_X ,n} ),
\]
\[
B_W = U_{B_W } \Lambda _{B_W } U_{B_W }^* ,\qquad\Lambda _{B_W }  = \mbox{diag}(\lambda _{B_W ,1} , \ldots ,\lambda _{B_W ,n} ).
\]

Let $\alpha$ be a given scaling factor. The parameter $\alpha$ regulates the strength of the embedding watermark into the host image. Now, embed the symmetric part of the watermark into the symmetric part of the host image by modifying the eigenvalues of $B_X$ and then obtain the symmetric part of the watermarked image $Y$ as detailed below:
\begin{equation}\label{hatt}
B_Y=U_{B_X } {\Lambda_Y}  U_{B_X }^*,\qquad{\Lambda_Y}=\Lambda _{B_X }+\alpha \Lambda _{B_W }.
\end{equation}

Hereinafter, reconstruction of the watermarked image $Y$ from its symmetric part is rendered. To do this, all the elements that are located above the main diagonal of the matrix $X$ are required. Since $Y$ can be construed as an acceptable approximation of $X$, let us consider all the elements above of the main diagonal of $Y$ as the elements of $X$. Obviously, the elements located below the main diagonal of the matrix $Y^T$ are also determined. In addition, it is concluded from the fact $2B_Y=Y+Y^T$ that the elements which are stated below the main diagonal of $Y$ can be obtained by subtracting the elements stated below the main diagonal of $2B_Y$ from the elements of $Y^T$ . See the above equation where $"\checkmark"$ and $"\times"$ denote the given and unknown entries, respectively.
\[
2B_Y  = \mathop {\left[ {\begin{array}{*{20}c}
  \times & {} & \checkmark  \\
   {} &  \ddots  & {}  \\
   \times & {} & \times  \\
\end{array}} \right]}\limits_Y  + \mathop {\left[ {\begin{array}{*{20}c}
   \times & {} & \times  \\
   {} &  \ddots  & {}  \\
   \checkmark & {} & \times  \\
\end{array}} \right]}\limits_{Y^T }.
\]

Finally, it is clear that the elements on the main diagonal of $Y$ are the main diagonal entries of $B_Y$.
It should be noted that by this procedure, the watermark is embedded only into the elements that are placed below the main diagonal of $X$, and its other elements are kept unchanged.

Furthermore, to reconstruct the watermarked image $Y$, we may reserve the elements located below the main diagonal of $X$ instead of the elements above the main diagonal and then follow a procedure similar to what is performed newly. Indeed, we can partition $X$ into several segments and reserve it provided that $Y$ can be reconstructed. This may be useful specially when unchanging some partitions of the host image is desirable during the digital image watermarking.

Now, as the next necessary step, we must be able to extract the watermark from the watermarked image. Our scheme for watermark extraction implicates obtaining the symmetric matrix $B_W$ and then reconstructing the watermark image $W$ from it. In order to compute the matrix $B_W$, the diagonal matrix $\Lambda_{B_X}$ and the unitary matrix $U_{B_W}$ are required. Keeping (\ref{hatt}) in mind, the following relations describe our proposed extraction technique for finding $B_W$
\[
\Lambda _W  = \frac{{ \Lambda_Y  - \Lambda _{B_X } }}{\alpha },\qquad B_W = U_{B_W } \Lambda _{B_W } U_{B_W }^*.
\]

Anew, reconstructing the $W$ needs the elements which are stated above (or below) the main diagonal of $W$ and can be performed similar to what was mentioned in the case of reconstructing the watermarked image $Y$. Let us once again remark that, if some segments of the watermark have more significance, then we can partition the watermark such that the important segments be preserved provided that reconstruction of the watermark is done.
% % % % % % % % % % % % % % % % % % % % % % % % % % % % % % % % % % % % % % % % % % % % % %
\subsection{Digital image watermarking method using both symmetric and skew-symmetric parts of $X$ and $W$}\label{sb2}

In the previous subsection, a digital image watermarking scheme which embedded a watermark image into almost half of a host image was introduced. That technique kept more than half of the elements of the host image unchanged, which might lead to some usages and advantages and, in particular, was a suitable method for small enough values of the scaling factor $\alpha$. However, embedding all of the watermark into half of the original image may cause the watermarked image to lose its reliability as compared to some other watermarking schemes for large enough scaling factors. To cope with this defect, we introduce a new method concerning both symmetric and skew-symmetric parts of the watermark and the host image in order to embed the watermark in the host image completely. This new watermarking technique provides for a remarkably high reliability.

As already mentioned, every matrix is equal to the sum of its symmetric and skew-symmetric parts. In the proposed technique, therefore, a given scaling factor is shared equally between the symmetric and skew-symmetric parts of the host image.  Keeping in mind the definitions of $B_X$ and $B_W$, $C_X$ and $C_W$ may be employed to denote the skew-symmetric parts of $X$ and $W$, respectively. The matrices $C_X$ and $C_W$ may be written as follows:
\[
C_X  = U_{C_X } \Lambda _{C_X } U_{C_X }^* ,\qquad\Lambda _{C_X }  = \mbox{diag}(\lambda _{C_X ,1} , \ldots ,\lambda _{C_X ,n} ),
\]
\[
C_W = U_{C_W } \Lambda _{C_W } U_{C_W }^* ,\qquad\Lambda _{C_W }  = \mbox{diag}(\lambda _{C_W ,1} , \ldots ,\lambda _{C_W ,n} ).
\]

Inspired by the method described in the previous subsection, for scaling factor $\alpha$, set
\[
\begin{gathered}
  \hat{C}_X  = U_{C_X } \hat \Lambda U_{C_X }^* ,\qquad\hat \Lambda  = \Lambda _{C_X }  + \frac{\alpha }
{2}\Lambda _{C_W},  \hfill \\
  \tilde{B}_X  = U_{B_X } \tilde \Lambda U_{B_X }^* ,\qquad\tilde \Lambda  = \Lambda _{B_X }  + \frac{\alpha }
{2}\Lambda _{B_W}.  \hfill \\ 
\end{gathered} 
\]
Therefore, the watermarked image $\mathcal{X}$ will be $\mathcal{X}=\tilde{B}_X+\hat{C}_X.$

Henceforth, extracting the watermark $W$ will be considered to be from the watermarked image $\mathcal{X}$. The extraction process needs two diagonal matrices, i.e., $\Lambda_{B_X}$ and $\Lambda_{C_X}$, and two unitary matrices,i.e., $U_{B_W}$ and $U_{C_W}$. In the presented watermark extraction scheme, at first, the matrices $\bar{\Lambda}_{B_W}$ and $\bar{\Lambda}_{C_W}$ are computed as the estimations of $\Lambda_{B_W}$ and $\Lambda_{C_W}$, respectively. Next, the extracted watermark $\mathcal{W}$ will be obtained by approximating the two matrices $B_W$ and $C_W$. Assume now, $\mathcal{X}$ is an acceptable approximation of the original image $X$. Form the symmetric and skew-symmetric parts of $\mathcal{X}$, namely $B_\mathcal{X}$ and $C_\mathcal{X}$, in the order already mentioned, the eigenvalue decompositions can be computed as shown below:
\[
\begin{gathered}
  B_\mathcal{X}=\dfrac{\mathcal{X}+\mathcal{X}^T}{2},\qquad B_\mathcal{X}= U_{B_\mathcal{X} } \Lambda _{B_\mathcal{X} } U_{B_\mathcal{X}}^*,  \hfill \\
C_\mathcal{X}=\dfrac{\mathcal{X}-\mathcal{X}^T}{2},\qquad C_\mathcal{X}= U_{C_\mathcal{X} } \Lambda _{C_\mathcal{X} } U_{C_\mathcal{X}}^*. \hfill \\ 
\end{gathered} 
\]

Now, the diagonal matrices $\bar{\Lambda}_{B_W}$ and $\bar{\Lambda}_{C_W}$ are obtained from the following relations:
\[\bar{\Lambda}_{B_W}=\dfrac{\Lambda _{B_\mathcal{X} }-\Lambda _{B_X }}{\alpha},\qquad\bar{\Lambda}_{C_W}=\dfrac{\Lambda _{C_\mathcal{X} }-\Lambda _{C_X }}{\alpha}.
\]
Also if we consider $\bar{B}_W$ and $\bar{C}_W$ estimations of $B_W$ and $C_W$, respectively, then
\[
\bar{B}_W = U_{B_W } \bar{\Lambda} _{B_W } U_{B_W }^* ,\qquad\mbox{and}\qquad\bar{C}_W = U_{C_W } \bar{\Lambda} _{C_W } U_{C_W }^*.
\]
Finally, the extracted watermark $\mathcal{W}$ will be $\mathcal{W}=\bar{B}_W +\bar{C}_W$.

We should note that none of the four required matrices in the extraction process provides information about the host image and the watermark image. This point improves the safety of the presented watermarking method. 

% % % % % % % % % % % % % % % % % % % % % % % % % % % % % % % %
\section{Experimental results}\label{exp}
In this section, the validity and the influence of the watermarking methods described in Section \ref{water} are examined by experiments. Moreover, let us note that the ideas presented in this paper, can be used immediately to establish a block-based watermarking technique. A block-based watermarking scheme concerns dividing both the watermark and the host image into non-overlapping blocks and applying a watermarking method to each block. The block-based SVD watermarking method is especially well known (e.g., see Section 4 of \cite{arathi}). 

In our experiments, the watermarking method explained in Subsection \ref{sb1} is considered for using either the symmetric or skew-symmetric parts of the host image and the watermark. Also, the watermarking scheme demonstrated in Subsection \ref{sb2} has been tested to prove its efficiency.

% % % % % % % % % % % % % % % % %

\begin{figure}
\centering
\includegraphics[width=0.7\linewidth]{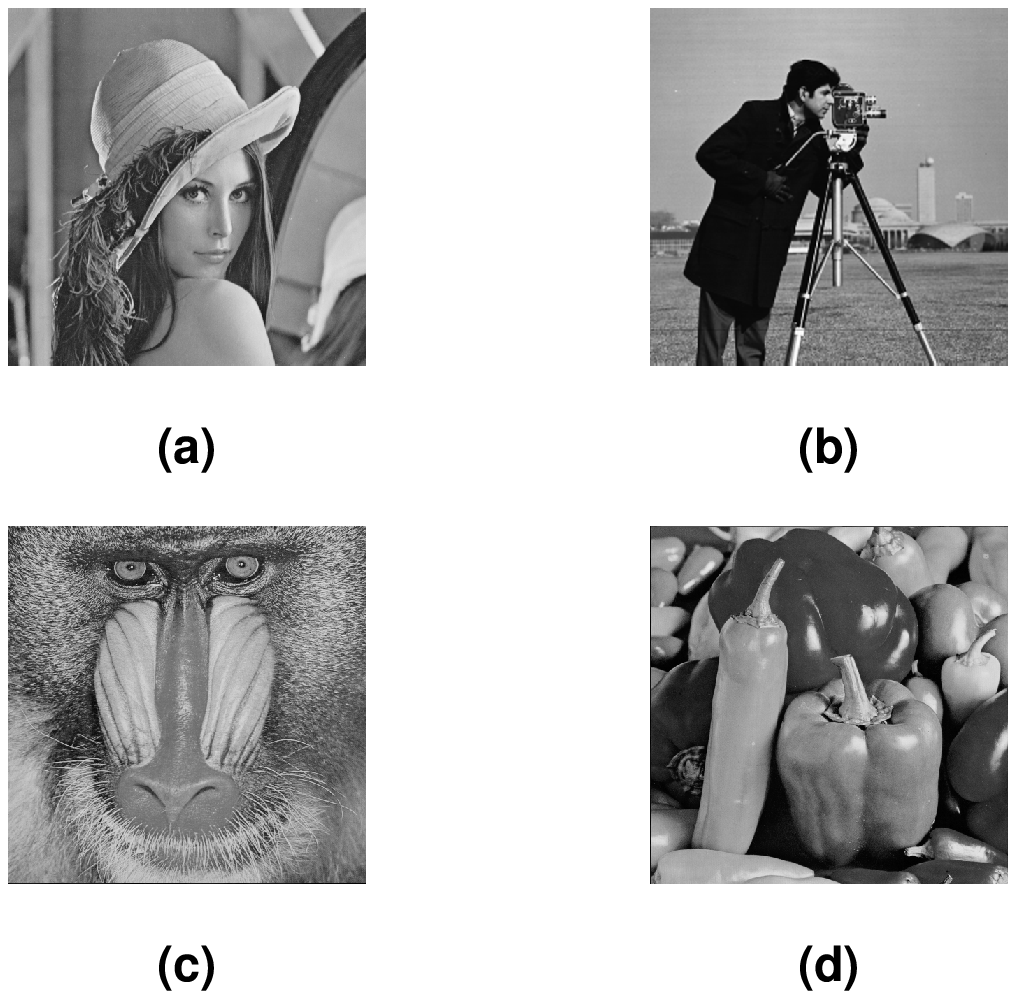}
\caption{The Host Images}
\label{fig:host_images}
\end{figure}

\begin{figure}
\centering
\includegraphics[width=0.7\linewidth]{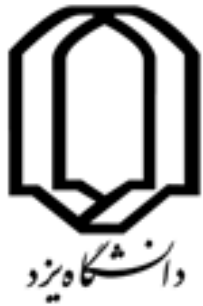}
\caption{The Watermark Image}
\label{fig:watermarks}
\end{figure}

The Peak Signal to Noise Ratio (PSNR) is calculated to measure the quality of the watermarked image. PSNR for gray scale images of size $M\times N$, whose pixels are represented with 8 bits, is computed by the following relations:
\[
PSNR = 10\log _{10}^{\frac{{255^2 }}{{MSE}}} ;\qquad MSE = \frac{1}{{MN}}\sum\limits_{i,j}^{} {\left| {X_{i,j}  - \mathcal{X}_{i,j} } \right|^2 },
\]
where $X_{i,j}$ and $\mathcal{X}_{i,j}$ denote the elements of the original and the watermarked images, respectively. In the above relation, MSE denotes the Mean Square Error between the host image and the watermarked image pixels. In the performed experiments, we considered four host $512 \times 512$ gray scale images (a) Lena, (b) Cameraman, (c) Baboon and (d) Peppers presented in Fig \ref{fig:host_images}. Also, our watermark image was Yazd University emblem of size $128\times 128$, as displayed in Fig \ref{fig:watermarks}.
The results are shown in Tables \ref{table:lena}, \ref{table:cameraman}, \ref{table:baboon} and \ref{table:peppers}, for some values of the scaling factor $\alpha$, for the host images Lena, Cameraman, Baboon and Peppers and the watermark image Yazd university emblem, respectively. In these tables, Method \#1 and Method \#2 stand for using the symmetric and the skew-symmetric parts of the host and watermark images, in the image watermarking schemes, respectively. In addition, Method \#3 denominates the watermarking technique described in Subsection \ref{sb2}. The results obtained by these techniques are compared to the results achieved by the digital image watermarking method using SVD which is denoted by Method \#4 in our tables.

\begin{figure}
\centering
\includegraphics[width=0.7\linewidth]{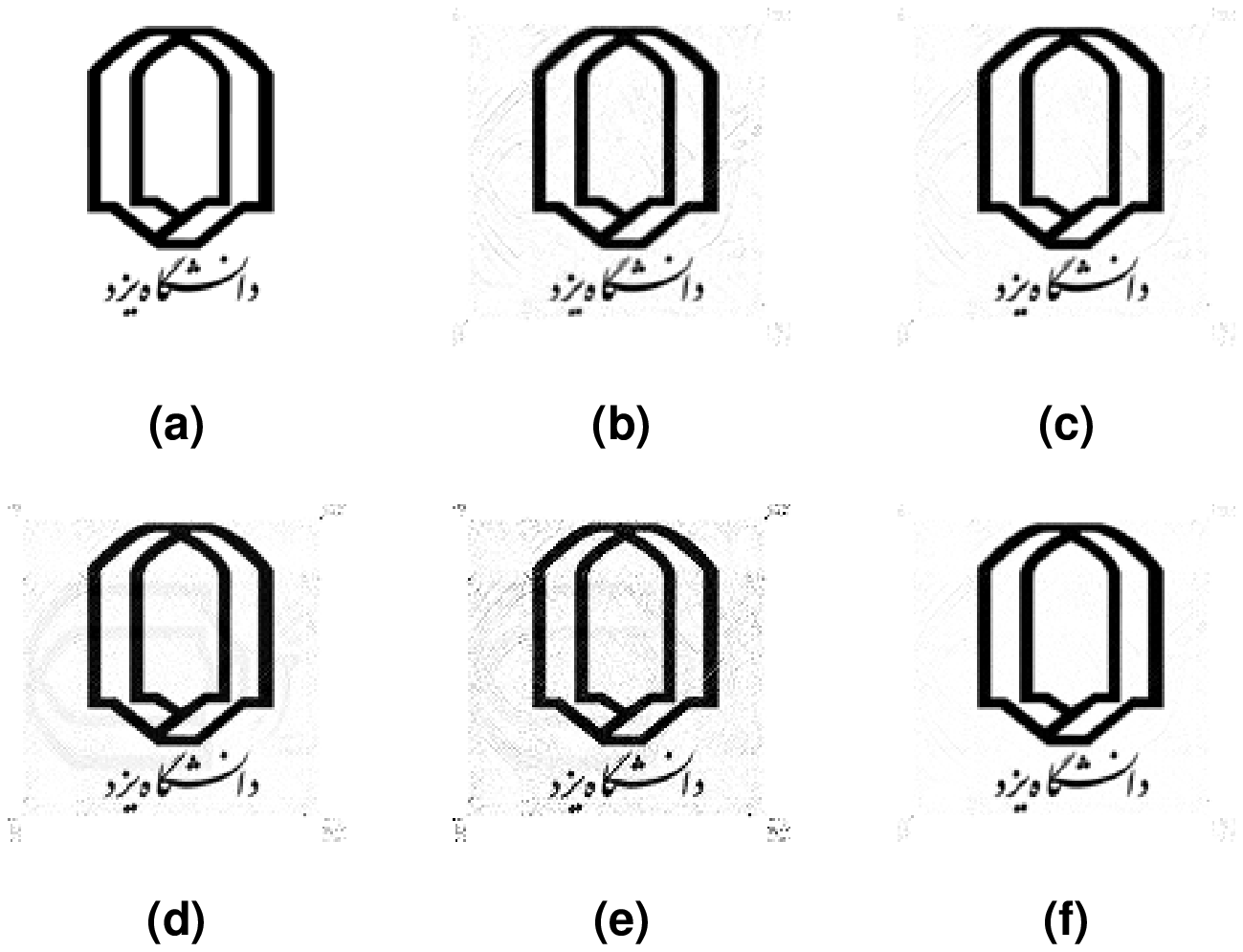}
\caption{Extracted watermark from watermarked image (a) Without any attack (b) After median filtering attack (c) After Gaussian low pas filtering attack (d) After JPEG compression (Quality=50) attack (e) After salt and pepper noise (density of 0.01) attack (f) After resizing (512-384-512) attack to the watermarked image. The used host image is Lena.}
\label{fig:attacked}
\end{figure}

Furthermore, the robustness of our watermarking schemes was evaluated in
the presence of different types of attacks. The Bit Error Rate (BER) was used to measure the robustness between the original watermark and the extracted watermark. BER is the percentage of bits with errors divided by the total number of bits that have been processed. We should note that the proposed methods described in Subsection \ref{sb1} are highly robust, and the BER results for those two methods almost zero under many types of attacks. Hence, the results corresponding to the other watermarking method, which was denoted by Method \#3, are only presented and compared with results obtained by Method \#4 to illustrate the robustness of our technique. The attacks used for examining the robustness of the method are median filtering, Gaussian low pass filtering, averaging filter, wiener filtering, JPEG compression, additive Gaussian white noise (AGWN), salt and
pepper noise, resizing and adjust image intensity values. Table \ref{table:bers} shows the BER results for the host image Lena and the watermark image Yazd University emblem under various attacks, where the window size is $3\times 3$. Moreover, Fig \ref{fig:attacked} shows the extracted watermark without any attack as well as after applying some sample attacks.

\begin{table}[ht]
\caption{PSNR results for Lena and Yazd University emblem} % title of Table
\centering % used for centering table
\begin{tabular}{c c c c c} % centered columns (5 columns)
\\ [-2 ex]
\hline\hline %inserts double horizontal lines
$\alpha$ & Method\#1 & Method\#2 & Method\#3 & Method \#4 \\ [.5ex] % inserts table
%heading
\hline % inserts single horizontal line
\\ [-2 ex]
    0.2000 & 72.1125 & 84.1944 & 74.8323 & 74.8500 \\
    0.4000 & 66.0919 & 78.1738 & 68.8117 & 68.8294 \\
    0.6000 & 62.5701 & 74.6520 & 65.2899 & 65.3076 \\ 
    0.8000 & 60.0713 & 72.1532 & 62.7911 & 62.8088 \\
    1.0000 & 58.1331 & 70.2150 & 60.8529 & 60.8706 \\
    1.2000 & 56.5495 & 68.6314 & 59.2693 & 59.2870 \\
    1.4000 & 55.2106 & 67.2925 & 57.9304 & 57.9481 \\
    1.6000 & 54.0507 & 66.1326 & 56.7705 & 56.7882 \\
    1.8000 & 53.0277 & 65.1096 & 55.7475 & 55.7652 \\
    2.0000 & 52.1125 & 64.1944 & 54.8323 & 54.8500 \\
[.5ex] % [1ex] adds vertical space
\hline %inserts single line
\end{tabular}
\label{table:lena} % is used to refer this table in the text
\end{table}
%rrrrrrrrrrrrrrrrrrrrrrrrrrrrrrrr
\begin{table}[ht]
\caption{PSNR results for Cameraman and  Yazd University emblem} % title of Table
\centering % used for centering table
\begin{tabular}{c c c c c} % centered columns (5 columns)
\\ [-2ex]
\hline\hline %inserts double horizontal lines
$\alpha$ & Method\#1 & Method\#2 & Method\#3 & Method \#4 \\ [.5ex] % inserts table
%heading
\hline % inserts single horizontal line
\\ [-2 ex]
    0.2000 & 72.1120 & 84.1944 & 74.8323 & 74.8500 \\
    0.4000 & 66.0914 & 78.1738 & 68.8117 & 68.8294 \\
    0.6000 & 62.5696 & 74.6520 & 65.2899 & 65.3076 \\
    0.8000 & 60.0708 & 72.1532 & 62.7911 & 62.8088 \\
    1.0000 & 58.1326 & 70.2150 & 60.8529 & 60.8706 \\
    1.2000 & 56.5490 & 68.6314 & 59.2693 & 59.2870 \\
    1.4000 & 55.2101 & 67.2925 & 57.9304 & 57.9481 \\
    1.6000 & 54.0502 & 66.1326 & 56.7705 & 56.7882 \\
    1.8000 & 53.0272 & 65.1096 & 55.7475 & 55.7652 \\
    2.0000 & 52.1120 & 64.1944 & 54.8323 & 54.8500 \\
[.5ex] % [1ex] adds vertical space
\hline %inserts single line
\end{tabular}
\label{table:cameraman} % is used to refer this table in the text
\end{table}
%rrrrrrrrrrrrrrrrrrrrrrrrrrrrrrrr
%rrrrrrrrrrrrrrrrrrrrrrrrrrrrrrrr
\begin{table}[ht]
\caption{PSNR results for Baboon and  Yazd University emblem} % title of Table
\centering % used for centering table
\begin{tabular}{c c c c c} % centered columns (5 columns)
\\ [-2ex]
\hline\hline %inserts double horizontal lines
$\alpha$ & Method\#1 & Method\#2 & Method\#3 & Method \#4 \\ [.5ex] % inserts table
%heading
\hline % inserts single horizontal line
\\ [-2 ex]
    0.2000 & 72.1121 & 84.1944 & 74.8323 & 74.8500 \\ 
    0.4000 & 66.0915 & 78.1738 & 68.8117 & 68.8294 \\
    0.6000 & 62.5696 & 74.6520 & 65.2899 & 65.3076 \\
    0.8000 & 60.0709 & 72.1532 & 62.7911 & 62.8088 \\
    1.0000 & 58.1327 & 70.2150 & 60.8529 & 60.8706 \\
    1.2000 & 56.5490 & 68.6314 & 59.2693 & 59.2870 \\
    1.4000 & 55.2101 & 67.2925 & 57.9304 & 57.9481 \\
    1.6000 & 54.0503 & 66.1326 & 56.7705 & 56.7882 \\
    1.8000 & 53.0272 & 65.1096 & 55.7475 & 55.7652 \\
    2.0000 & 52.1121 & 64.1944 & 54.8323 & 54.8500 \\
[.5ex] % [1ex] adds vertical space
\hline %inserts single line
\end{tabular}
\label{table:baboon} % is used to refer this table in the text
\end{table}
%rrrrrrrrrrrrrrrrrrrrrrrrrrrrrrrr
%rrrrrrrrrrrrrrrrrrrrrrrrrrrrrrrr
\begin{table}[ht]
\caption{PSNR results for Peppers and  Yazd University emblem} % title of Table
\centering % used for centering table
\begin{tabular}{c c c c c} % centered columns (5 columns)
\\ [-2ex]
\hline\hline %inserts double horizontal lines
$\alpha$ & Method\#1 & Method\#2 & Method\#3 & Method \#4 \\ [.5ex] % inserts table
%heading
\hline % inserts single horizontal line
\\ [-2 ex]

    0.2000 & 72.1127 & 84.1944 & 74.8323 & 74.8500 \\
    0.4000 & 66.0921 & 78.1738 & 68.8117 & 68.8294 \\
    0.6000 & 62.5703 & 74.6520 & 65.2899 & 65.3076 \\
    0.8000 & 60.0715 & 72.1532 & 62.7911 & 62.8088 \\
    1.0000 & 58.1333 & 70.2150 & 60.8529 & 60.8706 \\
    1.2000 & 56.5497 & 68.6314 & 59.2693 & 59.2870 \\
    1.4000 & 55.2108 & 67.2925 & 57.9304 & 57.9481 \\
    1.6000 & 54.0509 & 66.1326 & 56.7705 & 56.7882 \\
    1.8000 & 53.0279 & 65.1096 & 55.7475 & 55.7652 \\
    2.0000 & 52.1127 & 64.1944 & 54.8323 & 54.8500 \\
[.5ex] % [1ex] adds vertical space
\hline %inserts single line
\end{tabular}
\label{table:peppers} % is used to refer this table in the text
\end{table}
%rrrrrrrrrrrrrrrrrrrrrrrrrrrrrrrr
\begin{table}[ht]
\caption{BER results for the host image Lena and the watermark image Yazd University emblem under sample attacks} % title of Table
\centering
\begin{tabular}{l*{6}{c}r}
\\[-2ex]
Attack         &  Method \#3 & Method \#4 \\
\hline\hline
Median filtering & 0.2319 & 0.3479  \\
Gaussian low pas filtering              & 0.0793 & 0.2075 \\
Averaging filter    & 0.5249 & 1.0498 \\
Wiener filtering    & 0.2869 & 0.3296 \\
JPEG compression (Quality=70)     & 0.1587 & 0.0183 \\
JPEG compression (Quality=50)     & 0.2869  & 0.0305\\
AGWN    & 0.3540  & 0.8423\\
Salt and pepper noise (density of 0.01)    & 0.5005  & 1.2939\\
Resizing (512-384-512)  & 0.0427 & 0.1099  \\
Resizing (512-1024-512) & 0.0122 & 0.0427  \\
Adjust image intensity values    & 0.0427 & 66.0889 \\
\hline
\end{tabular}
\label{table:bers}
\end{table}
% % % % % % % % % % % % % % % % % % % % % % % % % % % % % % % % % %

% % % % % % % % % % % % % % %
\section{Conclusion}\label{conc}
In this paper, the eigenvalue decomposition of symmetric and skew-symmetric matrices, as two important kinds of normal matrices, and their properties were used to obtain digital image watermarking schemes. The proposed methods were straightforward and uncomplicated methods requiring clear and less computations compared to some exiting methods. This watermarking technique is the first one ever introduced and concerns the symmetric parts of a host image and a watermark, which can be also considered for the case of skew-symmetric parts of aforesaid images. It has proved to be highly robust method against many types of attacks. This method could keep half (or some selected segments) of the image unchanging and this feature may be applied for some usages. However, this watermarking method could embed the watermark into just half of the original image, which caused the watermarked image to lose its reliability for large enough values of the scaling factor $\alpha$. Hence, the second watermarking scheme using both symmetric and skew-symmetric parts of the original image and watermark was proposed. The experimental results illustrate the high reliability  and robustness of the methods. Finally, it is to be pointed out that these methods can be used to devise block-based watermarking techniques. Applying normal matrices to establish an image compression method will be the subject of our future research.

%%%%%%%%%%%%%%%%%%%%%%%%%%%%%%%%%%%%%%%%%%%%%%%%%%%%%%%%

\begin{thebibliography}{10}

{\small

\bibitem{arathi}
Ch. Arathi, A semi fragile image watermarking technique
using block based SVD, International Journal of Computer Science and Information Technologies. 3 (2012) 3644--3647.


\bibitem{normal1}
L. Elsner, Kh.D. Ikramov, Normal matrices: an update, Linear Algebra and its Applications, 285 (1998) 291--303.

%\bibitem{chandra}
%D.V.S. Chandra, Digital Image Watermarking Using Singular Value Decomposition, Proceedings of 45th IEEE
%Midwest Symposium on Circuits and Systems, 2002, pp. 264--267.


%\bibitem{chang}
%C.C. Chang, C.C. Lin, Y.S. Hu, An SVD oriented watermark embedding scheme with high qualities for the restored images, International Journal of Innovative Computing, Information and
%Control, 3 (2007)  pp. 609--620.


\bibitem{chang1}
C.C. Chang, P. Tsai, C.C. Lin, SVD-based digital image watermarking scheme, Pattern Recognition Letters,  26 (2005) 1577--1586.

\bibitem{chung}
K. Chung, W. Yang, Y. Huang, S. Wu, H. Yu-Chiao, On SVD-based watermarking algorithm, Applied Mathematics and Computation,  188 (2007) 54--57.


\bibitem{sym1}
J. Cullum, R.A. Willoughby, Computing eigenvalues of very large symmetric matrices An implementation of a Lanczos algorithm with no re-orthogonalization, Journal of Computational Physics, 44 (1981) 329--358.

\bibitem{sym2}
E.R. Davidson, The iterative calculation of a few of the lowest eigenvalues and corresponding eigenvectors of large real-symmetric matrices, Journal of Computational Physics, 17 (1975) 87--94.

\bibitem{sym4}
R.T. Gregory, Computing eigenvalues and eigenvectors of a symmetric matrix on the ILLIAC, Mathematical Tables and Other Aids to Computation, 7 (1953)  pp. 215--220.


\bibitem{normal2}
R. Grone, C.R. Johnson, E.M. Sa, H. Wolkowicz, Normal matrices, Linear Algebra and its Applications, 87 (1987) 213--225.

 
\bibitem{horn}
R.A. Horn, C.R. Johnson, Matrix Analysis, Cambridge University Press, Cambridge, 1985.

\bibitem{review}
D. Kannan, M. Gobi, An extensive research on robust digital image watermarking techniques: a review, International Journal of Signal and Imaging Systems Engineering, 8 (2015) 89--104.
%
%\bibitem{kum}
%P. Kumsawat, K. Attakitmongcol, A. Srikaew, An optimal robust digital image watermarking based on genetic algorithms in multiwavelet domain, Transactions on signal processing, 5 (2009) pp.
%42--51.

\bibitem{latif0}
A. Latif, An adaptive digital image watermarking scheme using fuzzy logic and tabu search, Journal of Information Hiding and Multimedia Signal Processing, 
4 (2013) 250--271.

\bibitem{latif}
A. Latif, A.R. Naghsh-Nilchi, Digital image watermarking based on parameters amelioration of parametric Slant-Hadamard transform using genetic algorithm, International Journal of Innovative Computing, Information and Control, 8 (2012) 1205--1220.


\bibitem{lin}
S. Lin, C. Chen, A robust DCT-based watermarking for copyright protection, IEEE Transactions on
consumer electronics, 46 (2000) pp. 415--421.


\bibitem{survey}
P. Parashar, R. Kumar Singh, A Survey: Digital Image Watermarking Techniques, International Journal of Signal Processing, Image Processing and Pattern Recognition, 7 (2014) pp. 111--124.

%\bibitem{potdar}
%V.M. Potdar, Song Han, E. Chang, A survey of digital image watermarking techniques," Industrial Informatics, 2005. INDIN '05. 2005 3rd IEEE International Conference on , vol., no., pp.709,716, 10-12 Aug. 2005

%\bibitem{pun}
%C. Pun, A novel DFT-based digital watermarking system for images, Int. Conf. on Signal Processing,l.2, 2006.

%\bibitem{normal3}
%H. Schneider, Theorems on normal matrices, Quarterly J. Math. Oxford Ser., 2 () 3241-249 (1952).

%\bibitem{sym3}
%I Shavitt, C.F. Bender, A. Pipano, R.P.  Hosteny, The iterative calculation of several of the lowest or highest eigenvalues and corresponding eigenvectors of very large symmetric matrices,\textit{ Journal of Computational Physics,}, 11 (1),1973, Pages 90-108.




}
\end{thebibliography}
\end{document}